# A DATA-DRIVEN BLOCK THRESHOLDING APPROACH TO WAVELET ESTIMATION

By T. Tony Cai[1] and Harrison H. Zhou[2]

*University of Pennsylvania and Yale University*

A data-driven block thresholding procedure for wavelet regression is proposed and its theoretical and numerical properties are investigated. The procedure empirically chooses the block size and threshold level at each resolution level by minimizing Stein's unbiased risk estimate. The estimator is sharp adaptive over a class of Besov bodies and achieves simultaneously within a small constant factor of the minimax risk over a wide collection of Besov Bodies including both the "dense" and "sparse" cases. The procedure is easy to implement. Numerical results show that it has superior finite sample performance in comparison to the other leading wavelet thresholding estimators.

**1. Introduction.** Consider the nonparametric regression model

$$y_i = f(t_i) + \sigma z_i, \qquad i = 1, 2, \ldots, n, \tag{1}$$

where $t_i = i/n$, $\sigma$ is the noise level and $z_i$'s are independent standard normal variables. The goal is to estimate the unknown regression function $f(\cdot)$ based on the sample $\{y_i\}$.

Wavelet methods have demonstrated considerable success in nonparametric regression. They achieve a high degree of adaptivity through thresholding of the empirical wavelet coefficients. Standard wavelet approaches threshold the empirical coefficients term by term based on their individual magnitudes. See, for example, Donoho and Johnstone (1994a), Gao (1998) and Antoniadis and Fan (2001). More recent work has demonstrated that block thresholding, which simultaneously keeps or kills all the coefficients in groups rather than individually, enjoys a number of advantages over the conventional term-by-term thresholding. Block thresholding increases estimation precision by

Received November 2006.

[1]Supported in part by NSF Grants DMS-03-06576 and DMS-06-04954.

[2]Supported in part by NSF Grant DMS-06-45676.

utilizing information about neighboring wavelet coefficients and allows the balance between variance and bias to be varied along the curve which results in adaptive smoothing. The degree of adaptivity, however, depends on the choice of block size and threshold level.

The idea of block thresholding can be traced back to Efromovich (1985) in orthogonal series estimators. In the context of wavelet estimation, global level-by-level thresholding was discussed in Donoho and Johnstone (1995) for regression and in Kerkyacharian, Picard and Tribouley (1996) for density estimation. But these block thresholding methods are not local, so they do not enjoy a high degree of spatial adaptivity. Hall, Kerkyacharian and Picard (1999) introduced a local blockwise hard thresholding procedure for density estimation with a block size of the order $(\log n)^2$ where $n$ is the sample size. Cai (1991) considered blockwise James–Stein rules and investigated the effect of block size and threshold level on adaptivity using an oracle inequality approach. In particular it was shown that a block size of order $\log n$ is optimal in the sense that it leads to an estimator which is both globally and locally adaptive. Cai and Silverman (2001) considered overlapping block thresholding estimators and Chicken and Cai (2005) applied block thresholding to density estimation.

The block size and threshold level play important roles in the performance of a block thresholding estimator. The local block thresholding methods mentioned above all have fixed block size and threshold and same thresholding rule is applied to all resolution levels regardless of the distribution of the wavelet coefficients. In the present paper, we propose a data-driven approach to empirically select both the block size and threshold at individual resolution levels. At each resolution level, the procedure, *SureBlock*, chooses the block size and threshold by minimizing Stein's Unbiased Risk Estimate (SURE). By empirically selecting both the block size and threshold and allowing them to vary from resolution level to resolution level, SureBlock has significant advantages over the more conventional wavelet thresholding estimators with fixed block sizes.

Both the numerical performance and asymptotic properties of SureBlock are studied in this paper. The SureBlock estimator is completely data-driven and easy to implement. A simulation study is carried out and the numerical results show that SureBlock has superior finite sample performance in comparison to the other leading wavelet estimators. More specifically, SureBlock uniformly outperforms both VisuShrink and SureShrink [Donoho and Johnstone (1994a, 1995)] in all 42 simulation cases in terms of the average squared error. SureBlock procedure is better than BlockJS [Cai (1991)] in 37 out of 42 cases.

The theoretical properties of SureBlock are considered in the Besov space formulation, that is, by now classical for the analysis of wavelet methods. Besov spaces, denoted by $B^{\alpha}_{p,q}$ and defined in Section 5, are a very rich class of function spaces which contain functions of inhomogeneous smoothness. The theoretical results show that SureBlock automatically adapts to the sparsity of the underlying wavelet coefficient sequence and enjoys excellent adaptivity over a wide range of Besov bodies. In particular, in the "dense case" $p \geq 2$ the SureBlock estimator is sharp adaptive over all Besov bodies $B_{p,q}^{\alpha}(M)$ with $p = q = 2$ and adaptively achieves within a factor of 1.25 of the minimax risk over Besov bodies $B_{p,q}^{\alpha}(M)$ for all $p \geq 2$, $q \geq 2$. At the same time the SureBlock estimator achieves simultaneously within a constant factor of the minimax risk over a wide collection of Besov bodies $B_{p,q}^{\alpha}(M)$ in the "sparse case" $p < 2$. These properties are not shared simultaneously by many commonly used fixed block size procedures such as VisuShrink [Donoho and Johnstone (1994a)], SureShrink [Donoho and Johnstone (1995)] or BlockJS [Cai (1991)].

The paper is organized as follows. In Section 2 we introduce the SureBlock method for the multivariate normal mean problem and derive oracle inequalities for the SureBlock estimator. The results developed in this section provide motivations and necessary technical tools for SureBlock in the wavelet regression setting. In Section 3, after a brief review of wavelets, the SureBlock procedure for the nonparametric regression is proposed. Section 4 discusses numerical implementation and compares the numerical performance of Sure-Block with those of VisuShrink [Donoho and Johnstone (1994a)], SureShrink [Donoho and Johnstone (1995)] and BlockJS [Cai (1991)]. Asymptotic properties of the SureBlock estimator are presented in Section 5. The proofs are given in Section 6.

**2. Estimation of a normal mean.** As mentioned in the Introduction, through an orthogonal discrete wavelet transform (DWT) the nonparametric regression problem can be turned into a problem of estimating the wavelet coefficients at individual resolution levels. The function estimation procedure as well as the analysis of the estimator become clear once the problem of estimating the wavelet coefficients at a given resolution level is well understood. In this section we shall treat the estimation problem at a single resolution level by considering a more generic problem, that of estimating the mean of a multivariate normal variable.

Suppose that we observe

$$(2) \qquad x_i = \theta_i + z_i, \qquad z_i \stackrel{\text{i.i.d.}}{\sim} N(0,1), \qquad i = 1, 2, \ldots, d,$$

and wish to estimate the mean vector $\theta = (\theta_1, \ldots, \theta_d)$ based on the observations $x = (x_1, \ldots, x_d)$ under the average mean squared error

$$(3) \qquad R(\widehat{\theta}, \theta) = d^{-1} \sum_{i=1}^{d} E(\widehat{\theta}_i - \theta_i)^2.$$

This normal mean problem occupies a central position in statistical estimation theory. Many methods have been introduced in the literature. In this section, with the application to wavelet function estimation in mind, we estimate the mean $\theta$ by a blockwise James–Stein estimator with block size $L$ and threshold level $\lambda$ chosen empirically by minimizing SURE. Oracle inequalities are developed in Section 2.2.

2.1. *SureBlock Procedure.* A block thresholding procedure thresholds the observations in groups and makes simultaneous decisions on all the means within a block. Let $L \geq 1$ be the possible length of each block, and $m = d/L$ be the number of blocks. (For simplicity we shall assume that $d$ is divisible by $L$ in the following discussion.) Fix a block size $L$ and a threshold level $\lambda$ and divide the observations $x_1, x_2, \ldots, x_d$ into blocks of size $L$. Let $\underline{x}_b = (x_{(b-1)L+1}, \ldots, x_{bL})$ represent observations in the $b$th block, and similarly $\underline{\theta}_b = (\theta_{(b-1)L+1}, \ldots, \theta_{bL})$ and $\underline{z}_b = (z_{(b-1)L+1}, \ldots, z_{bL})$. Let $S_b^2 = \|\underline{x}_b\|_2^2$ for $b = 1, 2, \ldots, m$. The blockwise James–Stein estimator is given by

$$(4) \qquad \widehat{\underline{\theta}}_b(\lambda, L) = \left(1 - \frac{\lambda}{S_b^2}\right)_+ \underline{x}_b, \qquad b = 1, 2, \ldots, m,$$

where $\lambda \geq 0$ is the threshold level. Block thresholding estimators depend on the choice of the block size $L$ and threshold level $\lambda$ which largely determines the performance of the resulting estimator. It is thus important to choose $L$ and $\lambda$ in an optimal way.

We shall select the block size $L$ and threshold level $\lambda$ by empirically minimizing SURE. Write $\widehat{\underline{\theta}}_b(\lambda, L) = \underline{x}_b + g(\underline{x}_b)$, where $g$ is a function from $\mathbb{R}^L$ to $\mathbb{R}^L$. Stein (1981) showed that when $g$ is weakly differentiable, then

$$E_{\underline{\theta}_b}\|\widehat{\underline{\theta}}_b(\lambda, L) - \underline{\theta}_b\|_2^2 = E_{\underline{\theta}_b}\{L + \|g\|_2^2 + 2\nabla \cdot g\}.$$

In our case, $g(\underline{x}_b) = (1 - \frac{\lambda}{S_b^2})_+ \underline{x}_b - \underline{x}_b$ is weakly differentiable. Simple calculations show $E_{\underline{\theta}_b}\|\widehat{\underline{\theta}}_b(\lambda, L) - \underline{\theta}_b\|_2^2 = E_{\underline{\theta}_b}(SURE(\underline{x}_b, \lambda, L))$, where

$$(5) \quad SURE(\underline{x}_b, \lambda, L) = L + \frac{\lambda^2 - 2\lambda(L-2)}{S_b^2} I(S_b^2 > \lambda) + (S_b^2 - 2L)I(S_b^2 \leq \lambda).$$

This implies that the total risk $E_\theta \|\widehat{\theta}(\lambda, L) - \theta\|_2^2 = E_\theta SURE(x, \lambda, L)$, where

$$(6) \qquad SURE(x, \lambda, L) = \sum_{b=1}^m SURE(\underline{x}_b, \lambda, L)$$

is an unbiased risk estimate. Our estimator is constructed through a hybrid method. Set $T_d = d^{-1}\sum(x_i^2 - 1)$, $\gamma_d = d^{-1/2}\log_2^{3/2} d$ and $\lambda^F = 2L \log d$. Let

$(\lambda^*, L^*)$ denote the minimizers of SURE with an additional restriction on the search range

$$(7) \qquad (\lambda^*, L^*) = \underset{\max\{L-2,0\} \leq \lambda \leq \lambda^F, 1 \leq L \leq d^{1/2}}{\arg\min} SURE(x, \lambda, L).$$

Define the estimator $\widehat{\theta}^*(x)$ of $\theta$ by

$$(8) \qquad \begin{aligned} \underline{\widehat{\theta}}^*{}_b &= \underline{\widehat{\theta}}_b(\lambda^*, L^*) \qquad \text{if } T_d > \gamma_d \quad \text{and} \\ \widehat{\theta}_i^* &= \left(1 - \frac{2\log d}{x_i^2}\right)_+ x_i \qquad \text{if } T_d \leq \gamma_d. \end{aligned}$$

We shall call this estimator the SureBlock estimator. When $T_d \leq \gamma_d$ the estimator is a degenerate block James–Stein estimator with block size $L = 1$. In this case the estimator is also called the nonnegative garrote estimator. See Breiman ([1995](#)) and Gao ([1998](#)). The SURE approach has also been used for the selection of the threshold level for fixed block size procedures, term-by-term thresholding ($L = 1$) in Donoho and Johnstone ([1995](#)) and block thresholding ($L = \log n$) in Chicken ([2005](#)).

REMARK. The hybrid scheme is used to guard against situations of extreme sparsity of the mean vector. See also Donoho and Johnstone ([1995](#)) and Johnstone ([1999](#)).

2.2. *Oracle inequalities.* We shall now consider the performance of the SureBlock estimator by comparing with that of ideal "estimators" equipped with an oracle. An oracle does not reveal the true estimand, but does "know" the optimal choice within a class of estimators. These ideal "estimators" are not true estimators in the statistical sense because the oracle depends on the unknown estimand. But the oracle risk of the "ideal estimators" provides a benchmark for the performance of estimators. It is desirable to have statistical estimators which can mimic the performance of the oracle. We shall consider two oracles: block thresholding oracle and linear shrinkage oracle. The oracle inequalities developed in this section are useful for showing the adaptivity results of SureBlock in the wavelet estimation setting. In particular, these results will be used in the proof of Theorem [3](#) given in Section [5](#).

*Block thresholding oracle.* Within the class of the block thresholding estimators, there is an "ideal estimator" which uses the optimal block size and threshold level so that the risk is minimized. The block thresholding oracle does not tell the true mean $\theta$, but "knows" the values of the ideal parameters,

$$(9) \qquad (\lambda^o, L^o) = \underset{0 \leq \lambda, 1 \leq L \leq d^{1/2}}{\arg\min} r(\lambda, L) = \underset{\max\{L-2,0\} \leq \lambda, 1 \leq L \leq d^{1/2}}{\arg\min} r(\lambda, L),$$

where $r(\lambda, L) = d^{-1} E \|\hat{\theta}(\lambda, L) - \theta\|_2^2$. Denote by $R_{block.oracle}(\theta)$ the oracle risk of the ideal block thresholding estimator $\hat{\theta}(\lambda^o, L^o)$, that is,

$$(10) \qquad R_{block.oracle}(\theta) = r(\lambda^o, L^o) = \inf_{\max\{L-2,0\} \leq \lambda, 1 \leq L \leq d^{1/2}} r(\lambda, L).$$

*Linear shrinkage oracle.* Linear shrinkage estimators have been commonly used in estimating a normal mean. A linear shrinker takes the form $\hat{\theta} = \gamma x$ where $0 \leq \gamma \leq 1$ is the shrinkage factor. The linear shrinkage oracle "knows" the ideal shrinkage factor $\gamma^*$ which equals $\|\theta\|_2^2/(\|\theta\|_2^2 + d)$. Simple calculations show that the risk of the ideal linear "estimator" $\tilde{\theta} = \gamma^* x$ is given by

$$R_{linear.oracle} = \frac{\|\theta\|_2^2}{\|\theta\|_2^2 + d}.$$

The following oracle inequalities show that the SureBlock estimator mimics the performance of both the block thresholding oracle and the linear shrinkage oracle.

THEOREM 1. *Let* $\{x_i, i = 1, \ldots, d\}$ *be given as in* (2) *and let* $\hat{\theta}^*$ *be the SureBlock estimator defined in* (8).

(a) (Block thresholding oracle.) *For some constant* $c > 0$,

$$(11) \qquad R(\hat{\theta}^*, \theta) \leq R_{block.oracle}(\theta) + cd^{-1/4}(\log d)^{5/2} \qquad \text{for all } \theta \in \mathbb{R}^d.$$

(b) (Linear shrinkage oracle.) *For some constant* $c > 0$,

$$R(\hat{\theta}^*, \theta) \leq R_{linear.oracle}(\theta) + cd^{-1/4}(\log d)^{5/2} \qquad \text{for all } \theta \in \mathbb{R}^d.$$

(c) *Set* $\mu_d = \|\theta\|_2^2/d$ *and* $\gamma_d = d^{-1/2} \log_2^{3/2} d$. *There exists some constant* $c > 0$ *such that for all* $\theta$ *satisfying* $\mu_d \leq \frac{1}{3}\gamma_d$

$$(12) \qquad R(\hat{\theta}^*, \theta) \leq d^{-1} \sum_i \theta_i^2 \wedge 2\log d + cd^{-1}(\log d)^{-1/2}.$$

Part (c) of Theorem 1 gives a risk bound of the SureBlock estimator in the case of $\theta$ being in a neighborhood of the origin. This bound is technically useful later for analysis in the wavelet function estimation setting. Parts (a) and (c) of Theorem 1 can be regarded as a generalization of Theorem 4 of Donoho and Johnstone (1995) from a fixed block size of one to variable block size. This generalization is important because it enables the resulting SureBlock estimator to be not only adaptively rate optimal over a wide collection of Besov bodies across both the dense $(p \geq 2)$ and sparse $(p < 2)$ cases, but also sharp adaptive over spaces where linear estimators can be asymptotically minimax. This property is not shared by fixed block size procedures such as VisuShrink, SureShrink and BlockJS, or the empirical Bayes estimator introduced in Johnstone and Silverman (2005).

In addition to the oracle inequalities given in Theorem 1, it is also interesting to consider the properties of the SureBlock estimator over $\ell_p$ balls,

$$\Theta_p(\tau) = \{\theta \in \mathbb{R}^d : \|\theta\|_p \le \tau\}. \tag{13}$$

THEOREM 2. Let $\{x_i, i = 1, \ldots, d\}$ be given as in (2) with $d \ge 4$ and let $\hat{\theta}^*$ be the SureBlock estimator defined in (8).

(a) (Adaptivity for dense signals.) For some constant $c > 0$ and for all $p \ge 2$

$$\sup_{\theta \in \Theta_p(\tau)} R(\hat{\theta}^*, \theta) \le \frac{\tau^2}{\tau^2 + d^{2/p}} + cd^{-1/4}(\log d)^{5/2}.$$

(b) (Adaptivity for moderate sparse signals.) For some constant $c > 0$ and for all $1 \le p \le 2$, $\sup_{\theta \in \Theta_p(\tau)} R(\hat{\theta}^*, \theta) \le cd^{-1}\tau^p(\log(d\tau^{-p}))^{1-p/2} + cd^{-1/4}(\log d)^{5/2}$.

(c) (Adaptive for a very sparse signals.) For $0 < p \le 2$ and $\tau < \frac{1}{\sqrt{3}}d^{1/4}\log_2^{3/4}d$, there is a constant $c > 0$ such that $\sup_{\theta \in \Theta_p(\tau)} R(\hat{\theta}^*, \theta) \le d^{-1}\tau^2 + cd^{-1} \times (\log d)^{-1/2}$.

## 3. The SureBlock procedure for wavelet regression.
Let $\{\phi, \psi\}$ be a pair of compactly supported father and mother wavelets with $\int \phi = 1$. Dilation and translation of $\phi$ and $\psi$ generate an orthonormal wavelet basis with an associated orthogonal DWT which transforms sampled data into the wavelet coefficient domain. A wavelet $\psi$ is called $r$-regular if $\psi$ has $r$ vanishing moments and $r$ continuous derivatives. See Daubechies (1992) and Strang (1992) for details on the DWT and compactly supported wavelets.

For simplicity in exposition, we work with periodized wavelet bases on $[0, 1]$. Let

$$\phi_{j,k}^p(x) = \sum_{l=-\infty}^{\infty} \phi_{j,k}(x-l), \qquad \psi_{j,k}^p(x) = \sum_{l=-\infty}^{\infty} \psi_{j,k}(x-l) \qquad \text{for } x \in [0, 1],$$

where $\phi_{j,k}(x) = 2^{j/2}\phi(2^j x - k)$ and $\psi_{j,k}(x) = 2^{j/2}\psi(2^j x - k)$. The collection $\{\phi_{j_0,k}^p, k = 1, \ldots, 2^{j_0}; \psi_{j,k}^p, j \ge j_0 \ge 0, k = 1, \ldots, 2^j\}$ is then an orthonormal basis of $L^2[0, 1]$, provided $j_0$ is large enough to ensure that the support of the wavelets at level $j_0$ is not the whole of $[0, 1]$. The superscript "$p$" will be suppressed from the notation for convenience. A square-integrable function $f$ on $[0, 1]$ can be expanded into a wavelet series

$$f(x) = \sum_{k=1}^{2^{j_0}} \xi_{j_0,k}\phi_{j_0,k}(x) + \sum_{j=j_0}^{\infty} \sum_{k=1}^{2^j} \theta_{j,k}\psi_{j,k}(x), \tag{14}$$

where $\xi_{j_0,k} = \langle f, \phi_{j_0,k} \rangle$ are the coefficients of the father wavelets at the coarsest level which represent the gross structure of the function $f$, and $\theta_{j,k} = \langle f, \psi_{j,k} \rangle$ are the wavelet coefficients which represent finer and finer structures as the resolution level $j$ increases.

Suppose we observe $Y = (y_1, \ldots, y_n)'$ as in (1) and suppose the sample size $n = 2^J$ for some integer $J > 0$. We use the standard device of the DWT to turn the function estimation problem into a problem of estimating wavelet coefficients. Let $\tilde{Y} = W \cdot n^{-1/2}Y$ be the DWTs of $n^{-1/2}Y$. Then $\tilde{Y}$ can be written as

$$(15) \quad \tilde{Y} = (\tilde{\xi}_{j_0,1}, \ldots, \tilde{\xi}_{j_0 2^{j_0}}, \tilde{y}_{j_0,1}, \ldots, \tilde{y}_{j_0,2^{j_0}}, \ldots, \tilde{y}_{J-1,1}, \ldots, \tilde{y}_{J-1,2^{J-1}})',$$

where $j_0$ is some fixed primary resolution level. Here $\tilde{\xi}_{j_0,k}$ are the gross structure terms, and $\tilde{y}_{j,k}$ are the empirical wavelet coefficients at level $j$ which represent fine structure at scale $2^j$. Since the DWT is an orthogonal transform, the $\tilde{y}_{j,k}$ are independent normal variables with standard deviation $\sigma_n = n^{-1/2}\sigma$. The mean of $\tilde{y}_{j,k}$, denoted by $\tilde{\theta}_{j,k}$, is the DWT of the sampled function $\{n^{-1/2}f(\frac{i}{n})\}$. Note that $\tilde{\theta}_{j,k}$ equals, approximately, the true wavelet coefficient $\theta_{j,k}$ of $f$. The approximation error is given in Lemma 4 in Section 6. Through the DWT, the nonparametric regression problem is then turned into a problem of estimating a high-dimensional normal mean vector.

3.1. *SureBlock for wavelet regression.* We now return to the nonparametric regression model (1). Denote by $\underline{\tilde{Y}}_j = \{\tilde{y}_{j,k} : k = 1, \ldots, 2^j\}$ and $\underline{\theta}_j = \{\theta_{j,k} : k = 1, \ldots, 2^j\}$ the empirical and true wavelet coefficients of the regression function $f$ at resolution level $j$. We apply the SureBlock procedure developed in Section 2 to the empirical wavelet coefficients $\underline{\tilde{Y}}_j$ level by level and then use the inverse DWT to obtain the estimate of the regression function. More specifically the SureBlock procedure for wavelet regression has the following steps.

1. Transform the data into the wavelet domain via DWT: $\tilde{Y} = W \cdot n^{-1/2}Y$.
2. At each resolution level $j$, estimate the wavelet coefficients using Sure-Block, that is,

$$(16) \quad \underline{\hat{\theta}}_j = \sigma_n \cdot \hat{\theta}^*(\sigma_n^{-1}\underline{\tilde{Y}}_j),$$

where $\sigma_n = n^{-1/2}\sigma$ and $\hat{\theta}^*$ is the SureBlock estimator given in (8). The estimate of the whole function $f$ is given by

$$(17) \quad \hat{f}^*(t) = \sum_{k=1}^{2^{j_0}} \tilde{\xi}_{j_0,k}\phi_{j_0,k}(t) + \sum_{j=j_0}^{J-1} \sum_{k=1}^{2^j} \hat{\theta}_{j,k}\psi_{j,k}(t).$$

3. The function at the sample points $\underline{f} = \{f(\frac{i}{n}) : i = 1, \ldots, n\}$ is estimated by the inverse transform of the denoised wavelet coefficients: $\underline{\hat{f}} = W^{-1} \cdot n^{1/2}\widehat{\Theta}$.

This procedure is easy to implement with good numerical performance. Theoretical results given in Section 5 show that $(L_j^*, \lambda_j^*)$ is optimal in the sense that the resulting estimator adaptively attains the exact minimax block thresholding risk asymptotically.

**4. Implementation and numerical results.** We now turn to the numerical performance of SureBlock. Proposition 1 below shows that for a given block size $L$ it suffices to search over the finite set $A$ for the threshold $\lambda$ which minimizes $SURE(x, \lambda, L)$. This makes the implementation of the SureBlock procedure easy. The result is also useful for the derivation of the theoretical results of SureBlock.

PROPOSITION 1. *Let $x_i, i = 1, \ldots, d$, and $SURE(x, \lambda, L)$ be given as in (2) and (6), respectively. Let the block size $L$ be given. Then the minimizer $\lambda$ of $SURE(x, \lambda, L)$ is an element of the set $A$ where $A \overset{\triangle}{=} \{x_i^2; 1 \leq i \leq d\} \cup \{0\}$ if $L = 1$, and $A \overset{\triangle}{=} \{S_i^2; S_i^2 \geq L - 2, 1 \leq i \leq m\} \cup \{L - 2\}$ if $L \geq 2$.*

The noise level $\sigma$ is assumed to be known in Section 2. In practice $\sigma$ needs to be estimated. As in Donoho and Johnstone (1994a) we estimate $\sigma$ based on the empirical coefficients at the highest resolution level by $\hat{\sigma} = \frac{1}{0.6745} \operatorname{median}(|n^{1/2}\tilde{y}_{J-1,k}| : 1 \leq k \leq 2^{J-1})$.

We now compare the numerical performance of SureBlock with that of VisuShrink [Donoho and Johnstone (1994a)], SureShrink [Donoho and Johnstone (1995)] and BlockJS [Cai (1991)]. VisuShrink thresholds empirical wavelet coefficients individually with a fixed threshold level. SureShrink is a soft thresholding procedure which selects the threshold at each resolution level by minimizing Stein's unbiased risk estimate. BlockJS is a block thresholding procedure with a fixed block size $\log n$ and a fixed threshold level. Each of these wavelet estimators has been shown to perform well numerically as well as theoretically. For further details see the original papers.

Six test functions, representing different levels of spatial variability, and various sample sizes, wavelets and signal to noise ratios are used for a systematic comparison of the four wavelet procedures. The test functions are plotted in the Appendix. Sample sizes ranging from $n = 256$ to $n = 16384$ and signal-to-noise ratios (SNR) from 3 to 7 were considered. The SNR is the ratio of the standard deviation of the function values to the standard deviation of the noise. Different combinations of wavelets and signal-to-noise ratios yield basically the same results. For reasons of space, we only report here the results for one particular case, using Daubechies' wavelet *Symmlet* 8 and SNR = 7. See Cai and Zhou (2005) for additional simulation results. We use the package WaveLab for simulations and the procedures MultiVisu (for VisuShrink) and MultiHybrid (for SureShrink) in WaveLab 802 are used (see http://www-stat.stanford.edu/~wavelab/).

Figure 1 reports the average squared errors (ASE) over 50 replications for the four thresholding estimators. SureBlock consistently outperforms both VisuShrink and SureShrink in all 42 simulation cases in terms of the ASE. SureBlock procedure is better than BlockJS about 88% of times (37 out of 42 cases). SureBlock fails to dominate BlockJS only for the test function "Doppler." For $n = 16384$ the risk ratio of SureBlock to BlockJS is $0.013/0.011 \approx 1.18$ and additional simulations show that the risk ratio goes to 1 as sample size increases. The main reason for BlockJS outperforming SureBlock in the case of "Doppler" is that at each resolution level the few significant wavelet coefficients all cluster together and this special structure greatly increases the accuracy of BlockJS. On the other hand, SureBlock is invariant to permutations of wavelet coefficients at any resolution level. Although SureBlock does not dominate BlockJS for "Doppler," the improvement of SureBlock over BlockJS is significant for other test functions. The simulation results show that, by empirically choosing the block size and threshold and allowing them to vary from resolution level to resolution level, the SureBlock estimator has significant numerical advantages over thresholding estimators with fixed block size $L = 1$ (VisuShrink or SureShrink) or $L = \log n$ (BlockJS). These numerical findings is consistent with the theoretical results given in Section 5.

Figure 2 shows an example of SureBlock applied to a noisy Bumps signal. The left panel is the noisy signal; the middle panel displays the empirical wavelet coefficients arranged according resolution levels; and the right panel is the SureBlock reconstruction (solid line) and the true signal (dotted line). In this example the block sizes chosen by SureBlock are 2, 3, 1, 5, 3, 5 and 1 from the resolution level $j = 3$ to level $j = 9$.

In addition to the comparison with other wavelet estimators, it is also instructive to compare the performance of SureBlock with the oracle risk $\frac{1}{n} \sum_{j,k} \theta_{j,k}^2 \wedge \sigma^2$, where $\theta_{j,k}$ are the true wavelet coefficients. Furthermore, to examine the advantage of empirically selecting block sizes, we compare the ASE of SureBlock with that of an estimator we call SureGarrote which empirically chooses the threshold at each level but fixes the block size $L = 1$. Figure 3 summarizes the numerical results for Doppler and Bumps with $n = 1024$, SNR ranging from 1 to 15 and 100 replications. SureBlock consistently outperforms SureGarrote in all cases. The ASE of SureGarrote is up to 40 percent higher than the corresponding ASE of SureBlock (see the right panels in Figure 3). Furthermore risk of SureBlock is within a small factor of the corresponding oracle risk. For Doppler the ratios of the ASE of

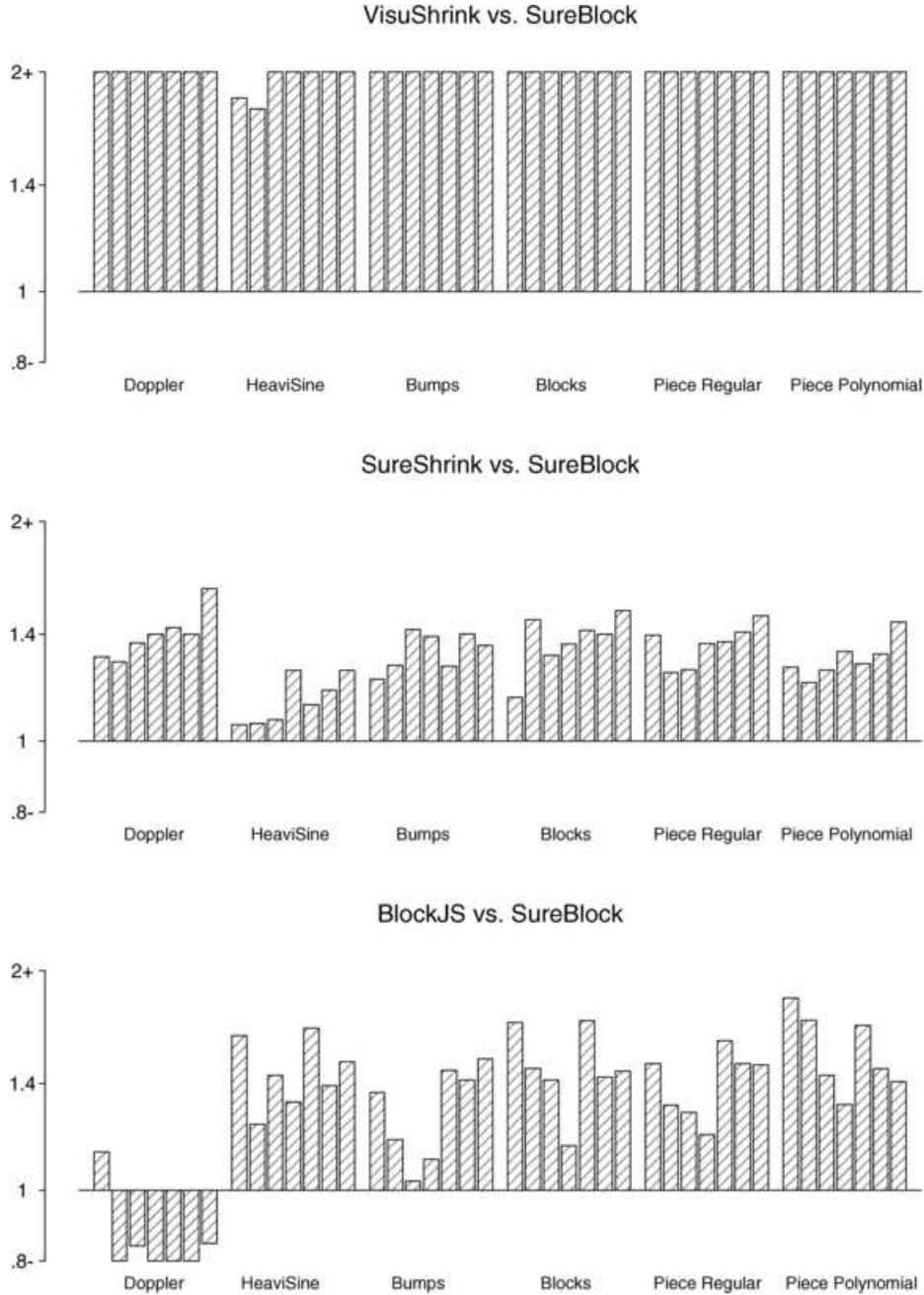

FIG. 1. *The vertical bars represent the ratios of the ASEs of estimators to the corresponding ASE of SureBlock. The higher the bar the better the relative performance of SureBlock. The bars are plotted on a log scale and are truncated at the value 2 of the original ratio. For each signal the bars are ordered from left to right by the sample sizes (n = 256 to 16382).*

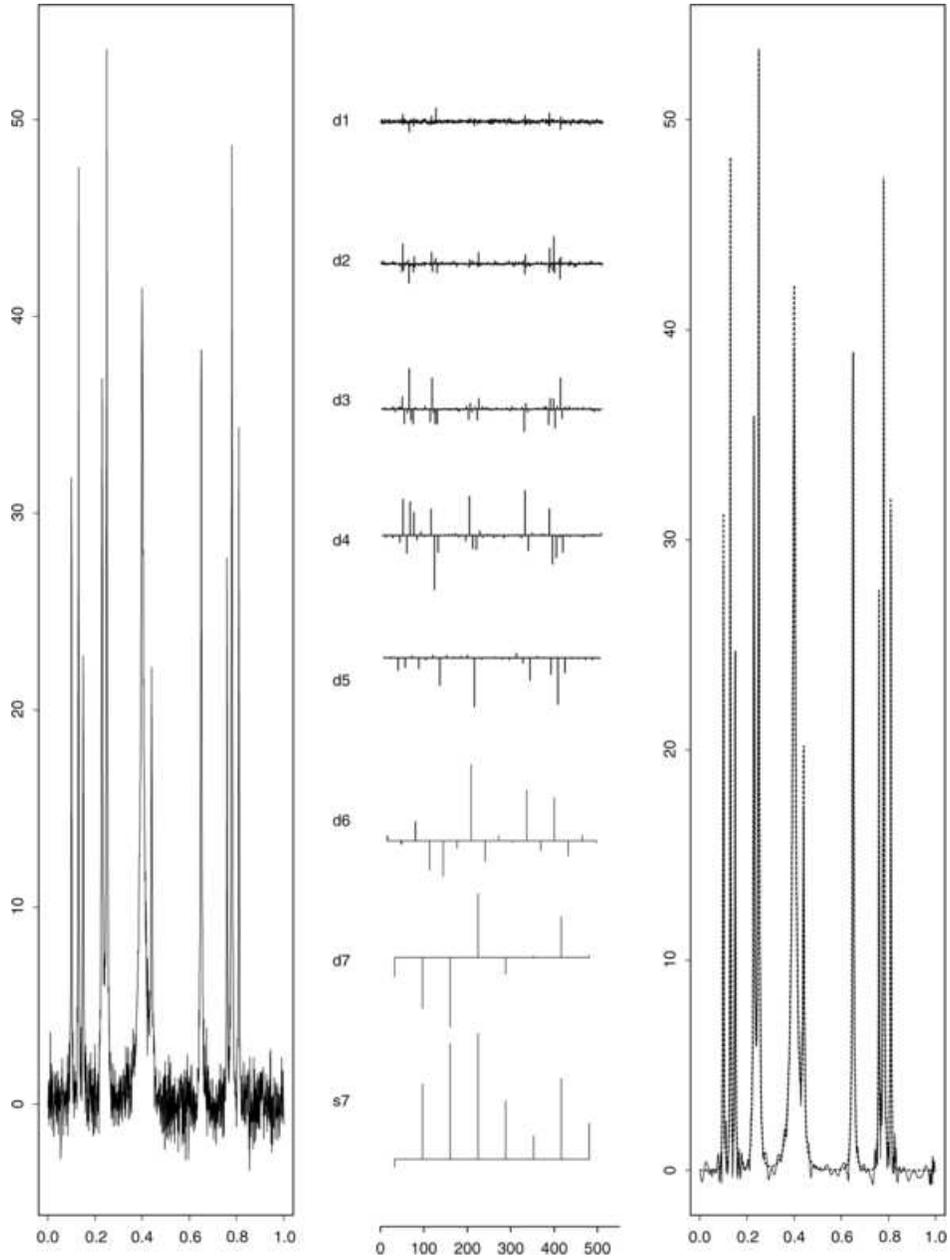

Fig. 2. *SureBlock procedure applied to a noisy Bumps signal.*

SureBlock and the oracle risk are between 2 to 2.7 and for Bumps the ratios are between 1.5 to 2.2 (see the left panels in Figure 3). In these simulations the block sizes chosen by SureBlock vary from 1 to 16, depending on the resolution levels. This example shows that the SureBlock procedure works well relative to the ideal oracle risk and empirically selecting block sizes improves the performance noticeably relative to the SureGarrote procedure. It would be interesting to carry out a more extensive numerical study to compare the performance of SureBlock with many other procedures including the empirical Bayes estimator of Johnstone and Silverman (2005). We leave this to future work.

**5. Theoretical properties of SureBlock.** We now turn to the theoretical properties of SureBlock for the nonparametric regression problem (1) under the integrated mean squared error $R(\hat{f}, f) = E\|\hat{f} - f\|_2^2$. The asymptotic results show that the SureBlock procedure is strongly adaptive.

Besov spaces are a very rich class of function spaces and contain as special cases many traditional smoothness spaces such as Hölder and Sobolev spaces. Roughly speaking, the Besov space $B_{p,q}^{\alpha}$ contains functions having $\alpha$ bounded derivatives in $L^p$ norm, the third parameter $q$ gives a finer gradation of smoothness. Full details of Besov spaces are given, for example, in Triebel (1983) and DeVore and Lorentz (1993). For a given $r$-regular mother wavelet $\psi$ with $r > \alpha$ and a fixed primary resolution level $j_0$, the Besov sequence norm $\|\cdot\|_{b_{p,q}^{\alpha}}$ of the wavelet coefficients of a function $f$ is then defined by

$$(18) \qquad \|f\|_{b_{p,q}^{\alpha}} = \|\underline{\xi}_{j_0}\|_p + \left( \sum_{j=j_0}^{\infty} (2^{js}\|\underline{\theta}_j\|_p)^q \right)^{1/q},$$

where $\underline{\xi}_{j_0}$ is the vector of the father wavelet coefficients at the primary resolution level $j_0$, $\underline{\theta}_j$ is the vector of the wavelet coefficients at level $j$, and $s = \alpha + \frac{1}{2} - \frac{1}{p} > 0$. Note that the Besov function norm of index $(\alpha, p, q)$ of a function $f$ is equivalent to the sequence norm (18) of the wavelet coefficients of the function. See Meyer (1992). The Besov body $B_{p,q}^{\alpha}(M)$ is defined by $B_{p,q}^{\alpha}(M) = \{f : \|f\|_{b_{p,q}^{\alpha}} \leq M\}$. The minimax risk of estimating $f$ over the Besov body $B_{p,q}^{\alpha}(M)$ is

$$(19) \qquad R^*(B_{p,q}^{\alpha}(M)) = \inf_{\hat{f}} \sup_{f \in B_{p,q}^{\alpha}(M)} E\|\hat{f} - f\|_2^2.$$

Donoho and Johnstone (1998) show that the minimax risk $R^*(B_{p,q}^{\alpha}(M))$ converges to 0 at the rate of $n^{-2\alpha/(1+2\alpha)}$ as $n \to \infty$.

The blockwise James–Stein estimation of the wavelet coefficients and the corresponding function $f$ is determined by the block size $L_j$ and threshold level $\lambda_j$ of each resolution $j$. Let $\underline{L} = (L_j)_{j \geq j_0}$ with $1 \leq L_j \leq 2^{j/2}$, and $\underline{\lambda} = (\lambda_j)_{j \geq j_0}$ with $\lambda_j \geq 0$. Let $\hat{f}_{\underline{L},\underline{\lambda}}$ be the corresponding estimator of $f$. The minimax risk among all block James–Stein estimators with all possible block sizes $\underline{L}$ and threshold levels $\underline{\lambda}$ is

$$(20) \qquad R_T^*(B_{p,q}^\alpha(M)) = \inf_{\widehat{f}_{\underline{L},\underline{\lambda}}} \sup_{f \in B_{p,q}^\alpha(M)} E\|\widehat{f}_{\underline{L},\underline{\lambda}} - f\|_2^2$$

and equivalently

$$R_T^*(B_{p,q}^\alpha(M)) = \inf_{\lambda_j \geq 0, 1 \leq L_j \leq 2^{j/2}} \sup_{\theta \in B_{p,q}^\alpha(M)} E \sum_{j=j_0}^\infty \|\widehat{\underline{\theta}}_j(\lambda_j, L_j) - \underline{\theta}_j\|_2^2.$$

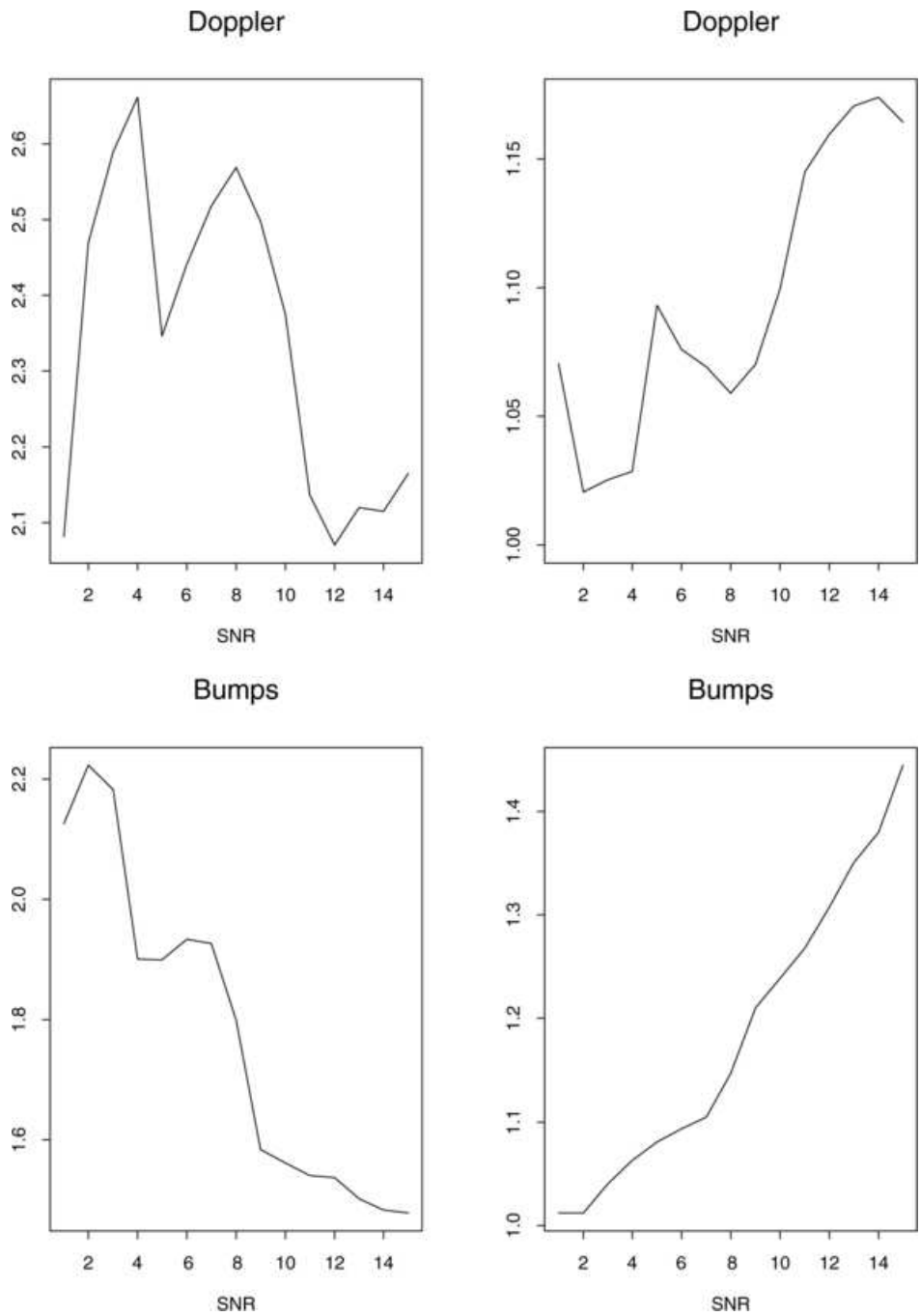

Fig. 3.  *Left panels: the ratios of the ASE of SureBlock and the oracle risk. Right panels: the ratios of the ASE of SureGarrote and the ASE of SureBlock.*

We shall call $R_T^*(B_{p,q}^\alpha(M))$ the minimax block thresholding risk. It is clear that $R_T^*(B_{p,q}^\alpha(M)) \geq R^*(B_{p,q}^\alpha(M))$. Theorems 4 and 5 below show $R_T^*(B_{p,q}^\alpha(M))$ is within a small constant factor of the minimax risk $R^*(B_{p,q}^\alpha(M))$. The following theorem shows that SureBlock adaptively attains the exact minimax block thresholding risk $R_T^*(B_{p,q}^\alpha(M))$ asymptotically over a wide range of Besov bodies.

THEOREM 3. *Suppose the mother wavelet $\psi$ is $r$-regular. Let $\widehat{f}^*$ be the SureBlock estimator of $f$ defined in (17). Then*

$$(21) \qquad \sup_{f \in B_{p,q}^\alpha(M)} E_f \|\widehat{f}^* - f\|_2^2 \leq R_T^*(B_{p,q}^\alpha(M))(1 + o(1))$$

*for $1 \leq p, q \leq \infty$, $0 < M < \infty$, and $r \geq \alpha > 4(\frac{1}{p} - \frac{1}{2})_+ + \frac{1}{2}$ with $\frac{2\alpha^2 - 1/6}{1 + 2\alpha} > 1/p$.*

Theorem 3 is proved in Section 6. The main technical tools for the proof are the oracle inequalities for SureBlock developed in Section 2.2.

Theorems 4 and 5 below make it clear that the SureBlock procedure is indeed nearly optimally adaptive over a wide collection of Besov bodies $B_{p,q}^\alpha(M)$ including both the dense $(p \geq 2)$ and sparse $(p < 2)$ cases. The estimator is asymptotically sharp adaptive over Besov bodies with $p = q = 2$ in the sense that it adaptively attains both the optimal rate and optimal constant. Over Besov bodies with $p \geq 2$ and $q \geq 2$ SureBlock adaptively achieves within a factor of 1.25 of the minimax risk. At the same time the maximum risk of the estimator is simultaneously within a constant factor of the minimax risk over a collection of Besov bodies $B_{p,q}^\alpha(M)$ in the sparse case of $p < 2$.

THEOREM 4. *Suppose $\psi$ is $r$-regular.* (i) *SureBlock is adaptively sharp minimax over Besov bodies $B_{2,2}^\alpha(M)$ for all $M > 0$ and $r \geq \alpha > 0.88$, that is,*

$$(22) \qquad \sup_{f \in B_{2,2}^\alpha(M)} E_f \|\widehat{f}^* - f\|_2^2 \leq R^*(B_{2,2}^\alpha(M))(1 + o(1)).$$

(ii) *SureBlock is adaptively, asymptotically within a factor of* 1.25 *of the minimax risk over Besov bodies* $B_{p,q}^{\alpha}(M)$,

$$(23) \qquad \sup_{f \in B_{p,q}^{\alpha}(M)} E_f \|\widehat{f}^* - f\|_2^2 \le 1.25 R^*(B_{p,q}^{\alpha}(M))(1 + o(1))$$

*for all* $p \ge 2$, $q \ge 2$, $M > 0$ *and* $\frac{2\alpha^2 - 1/6}{1 + 2\alpha} > 1/p$ *with* $r \ge \alpha > 1/2$.

For the sparse case $p < 2$, the SureBlock estimator is also simultaneously within a small constant factor of the minimax risk.

THEOREM 5. *Suppose* $\psi$ *is* $r$-regular. SureBlock *is asymptotically minimax up to a constant factor* $G(p \wedge q)$ *over a large range of Besov bodies with* $1 \le p, q \le \infty$, $0 < M < \infty$, *and* $r \ge \alpha > 4(\frac{1}{p} - \frac{1}{2})_+ + \frac{1}{2}$ *with* $\frac{2\alpha^2 - 1/6}{1 + 2\alpha} > 1/p$. *That is,*

$$(24) \qquad \sup_{f \in B_{p,q}^{\alpha}(M)} E_f \|\widehat{f}^* - f\|_2^2 \le G(p \wedge q) \cdot R^*(B_{p,q}^{\alpha}(M))(1 + o(1)),$$

*where* $G(p \wedge q)$ *is a constant depending only on* $p \wedge q$.

**6. Proofs.** Throughout this section, without loss of generality, we shall assume the noise level $\sigma = 1$. We first prove Theorem 1 and then use it as the main tool to prove Theorem 3. The proofs of Theorems 4 and 5 and Proposition 1 are given later.

6.1. *Notation and preparatory results.* Before proving the main theorems, we need to introduce some notation and collect a few technical results. The proofs of some of these preparatory results are long. For reasons of space these proofs are omitted here. We refer interested readers to Cai and Zhou (2005) for the complete proofs.

Consider the normal mean problem (2) with $\sigma = 1$. For a given block size $L$ and threshold level $\lambda$, set $r_b(\lambda, L) = E_{\underline{\theta}_b} \|\widehat{\underline{\theta}}_b(\lambda, L) - \underline{\theta}_b\|^2$ and define $r(\lambda, L) = \frac{1}{d} \sum_{b=1}^m r_b(\lambda, L) = ED(\lambda, L)$, where $D(\lambda, L) = \frac{1}{d} \sum_{b=1}^m \|\widehat{\underline{\theta}}_b(\lambda, L) - \underline{\theta}_b\|_2^2$. Set

$$(25) \qquad \widetilde{R}(\theta) = \inf_{\lambda \le \lambda^F, 1 \le L \le d^{1/2}} r(\lambda, L) = \inf_{\max\{L-2,0\} \le \lambda \le \lambda^F, 1 \le L \le d^{1/2}} r(\lambda, L).$$

The difference between $\widetilde{R}(\theta)$ and $R_{block.oracle}(\theta)$ defined in (10) is that the search range for the threshold $\lambda$ in $\widetilde{R}(\theta)$ is restricted to be at most $\lambda^F$. The result given below shows that the effect of this restriction is negligible for any block size $L$.

LEMMA 1. *For any fixed* $\eta > 0$, *there exists a constant* $C_\eta > 0$ *such that for all* $\theta \in R^d$,

$$\widetilde{R}(\theta) - R_{block.oracle}(\theta) \le C_\eta d^{\eta - 1/2}.$$

The following lemma is adapted from Donoho and Johnstone ([1995](#)) and is used in the proof of Theorem [1](#).

LEMMA 2.   *Let* $T_d = d^{-1} \sum (x_i^2 - 1)$ *and* $\mu_d = d^{-1} \|\theta\|_2^2$. *If* $\gamma_d^2 d / \log d \to \infty$, *then*

$$\sup_{\mu_d \geq 3\gamma_d} (1 + \mu_d) P(T_d \leq \gamma_d) = o(d^{-1/2}).$$

We also need the following bounds for the loss of the SureBlock estimator. This bound is used in the proof of Theorem [3](#).

LEMMA 3.   *Let* $\{x_i : i = 1, \ldots, d\}$ *be given as in* ([2](#)). *Then*

$$\|\hat{\theta}^* - \theta\|_2^2 \leq 4d \log d + 2\|z\|_2^2. \tag{26}$$

Finally we develop a key technical result for the proof of Theorem [1](#). Set

$$U(\lambda, L) \triangleq \frac{1}{d} SURE(x, \lambda, L)$$

$$= 1 + \frac{1}{d} \sum_{b=1}^{m} \left( \frac{\lambda^2 - 2\lambda(L-2)}{S_b^2} I(S_b^2 > \lambda) + (S_b^2 - 2L) I(S_b^2 \leq \lambda) \right).$$

Note that both $D(\lambda, L)$ and $U(\lambda, L)$ have expectation $r(\lambda, L)$.

The goal is to show that the minimizer $(\lambda^S, L^S)$ of $U(\lambda, L)$ is asymptotically the ideal threshold level and block size. The key step is to show that $\Delta_d = |ED(\lambda^S, L^S) - \inf_{\lambda, L} r(\lambda, L)|$ is negligible for $\max\{L - 2, 0\} \leq \lambda \leq \lambda^F$ and $1 \leq L \leq d^{1/2}$. Note that for two functions $g$ and $h$ defined on the same domain, $|\inf_x g(x) - \inf_x h(x)| \leq \sup_x |g(x) - h(x)|$. Hence, $|U(\lambda^S, L^S) - \inf_{\lambda, L} r(\lambda, L)| = |\inf_{\lambda, L} U(\lambda, L) - \inf_{\lambda, L} r(\lambda, L)| \leq \sup_{\lambda, L} |U(\lambda, L) - r(\lambda, L)|$ and consequently

$$\Delta_d \leq E \bigg| D(\lambda^S, L^S) - r(\lambda^S, L^S) + r(\lambda^S, L^S)$$

$$- U(\lambda^S, L^S) + U(\lambda^S, L^S) - \inf_{\lambda, L} r(\lambda, L) \bigg| \tag{27}$$

$$\leq E \sup_{\lambda, L} |D(\lambda, L) - r(\lambda, L)| + 2E \sup_{\lambda, L} |r(\lambda, L) - U(\lambda, L)|.$$

The upper bounds for the two terms on the RHS of ([27](#)) is given as follows.

PROPOSITION 2.   *Let* $\lambda^F = 2L \log d$. *Uniformly in* $\theta \in \mathbb{R}^d$, *we have*

$$E_\theta \sup_{\max\{L-2, 0\} \leq \lambda \leq \lambda^F, 1 \leq L \leq d^{1/2}} |U(\lambda, L) - r(\lambda, L)| \leq cd^{-1/4} (\log d)^{5/2}, \tag{28}$$

$$E_\theta \sup_{\max\{L-2, 0\} \leq \lambda \leq \lambda^F, 1 \leq L \leq d^{1/2}} |D(\lambda, L) - r(\lambda, L)| \leq cd^{-1/4} (\log d)^{5/2}. \tag{29}$$

The following result, which is crucial for the proof of Theorem 5, plays the role similar to that of Proposition 13 in Donoho and Johnstone (1994b).

PROPOSITION 3. *Let* $X \sim N(\mu, 1)$ *and let* $\mathcal{F}_p(\eta)$ *denote the probability measures* $F(d\mu)$ *satisfying* $\int |\mu|^p F(d\mu) \leq \eta^p$. *Let* $\overline{r}(\delta_\lambda^g, \eta) = \sup_{\mathcal{F}_p(\eta)} \{ E_F r_g(\mu) : \int |\mu|^p F(d\mu) \leq \eta^p \}$ *where* $r_g(\mu) = E_\mu (\delta_\lambda^g(x) - \mu)^2$ *and* $\delta_\lambda^g(x) = (1 - \frac{\lambda^2}{x^2})_+ x$. *Let* $p \in (0, 2)$ *and* $\lambda = \sqrt{2 \log \eta^{-p}}$, *then* $\overline{r}(\delta_\lambda^g, \eta) \leq 2\eta^p \lambda^{2-p} (1 + o(1))$ *as* $\eta \to 0$.

The following lemma bounds the approximation errors between the mean of the empirical wavelet coefficient and the true wavelet coefficient of $f \in B_{p,q}^\alpha(M)$. Set $\beta = \alpha - 1/p$, which is positive under the assumption $\frac{2\alpha^2 - 1/6}{1 + 2\alpha} > 1/p$ in Theorems 3, 4 and 5.

LEMMA 4. *Let* $\tilde{\theta} = (\tilde{\theta}_{j,k})$ *be the DWT of the sampled function* $\{n^{-1/2} f(\frac{k}{n})\}$ *with* $n = 2^J$ *and let* $\theta_{j,k} = \int f(x) \psi_{j,k}(x) \, dx$. *Then* $\|\tilde{\theta} - \theta\|_2^2 \leq C n^{-2\beta}$.

REMARK 1. Lemma 4 implies

$$
\inf_{\lambda_j \geq 0, 1 \leq L_j \leq 2^{j/2}} \sup_{\theta \in B_{p,q}^\alpha(M)} E \sum_{j=j_0}^\infty \|\underline{\hat{\theta}}_j(\lambda_j, L_j) - \underline{\tilde{\theta}}_j\|_2^2
$$
$$
= (1 + o(1)) R_T^*(B_{p,q}^\alpha(M))
\tag{30}
$$

and for $p \geq 2$ and $q > 2$

$$
\sup_{\theta \in B_{p,q}^\alpha(M)} \sum_{j=j_0}^\infty \sum_k \left( \frac{\theta_{j,k}^2/n}{\theta_{j,k}^2 + 1/n} \right)
$$
$$
= (1 + o(1)) \sup_{\theta \in B_{p,q}^\alpha(M)} \sum_{j=j_0}^\infty \sum_k \left( \frac{\tilde{\theta}_{j,k}^2/n}{\tilde{\theta}_{j,k}^2 + 1/n} \right)
\tag{31}
$$

under the assumption $\frac{2\alpha^2 - 1/6}{1 + 2\alpha} > 1/p$ which implies $\beta > \alpha/(2\alpha + 1)$. The argument for (30) is as follows. Write

$$
E \sum_{j=j_0}^\infty \|\underline{\hat{\theta}}_j(\lambda_j, L_j) - \underline{\tilde{\theta}}_j\|_2^2 - E \sum_{j=j_0}^\infty \|\underline{\hat{\theta}}_j(\lambda_j, L_j) - \underline{\theta}_j\|_2^2
$$
$$
= \|\tilde{\theta} - \theta\|_2^2 + 2E \sum_{j=j_0}^\infty \langle \underline{\hat{\theta}}_j(\lambda_j, L_j) - \underline{\theta}_j, \underline{\theta}_j - \underline{\tilde{\theta}}_j \rangle.
$$

From Lemma 4, $\|\tilde{\theta} - \theta\|_2^2 \leq C n^{-2\beta} = o(n^{-2\alpha/(2\alpha+1)})$. The Cauchy–Schwarz inequality implies

$$
\sup_{\theta \in B_{p,q}^\alpha(M)} E \sum_{j=j_0}^\infty \langle \underline{\hat{\theta}}_j(\lambda_j, L_j) - \underline{\tilde{\theta}}_j, \underline{\theta}_j - \underline{\tilde{\theta}}_j \rangle
$$

$$= \sup_{\theta \in B_{p,q}^\alpha(M)} \|\tilde{\theta} - \theta\|_2 \sqrt{E \sum_{j=j_0}^{\infty} \|\widehat{\underline{\theta}}_j(\lambda_j, L_j) - \underline{\theta}_j\|_2^2},$$

which is $o(n^{-2\alpha/(2\alpha+1)})$, since $\|\tilde{\theta} - \theta\|_2 = O(n^{-\beta})$ with $\beta > \alpha/(2\alpha + 1)$ and $R_T^*(B_{p,q}^\alpha(M)) \leq Cn^{-2\alpha/(1+2\alpha)} \log n$ from Cai (1991) in which $L_j = \log n$ and $\lambda_j = 4.505 \log n$. We know $R_T^*(B_{p,q}^\alpha(M)) \geq R^*(B_{p,q}^\alpha(M)) \geq Cn^{-2\alpha/(1+2\alpha)}$ from Donoho and Johnstone (1998). Thus (30) is established. The argument for (31) is similar.

In the following proofs we will denote by $C$ a generic constant that may vary from place to place.

6.2. *Proof of Theorem 1.* The proof of Theorem 1 is similar to but more involved than those of Theorem 4 of Donoho and Johnstone (1995) and Theorem 2 of Johnstone (1999) because of variable block size.

Set $T_d = d^{-1} \sum(x_i^2 - 1)$ and $\gamma_d = d^{-1/2} \log_2^{3/2} d$. Define the event $A_d = \{T_d \leq \gamma_d\}$ and decompose the risk of SureBlock into two parts:

$$R(\widehat{\theta}^*, \theta) = d^{-1} E_\theta \{\|\widehat{\theta}^* - \theta\|_2^2 I(A_d)\} + d^{-1} E_\theta \{\|\widehat{\theta}^* - \theta\|_2^2 I(A_d^c)\}$$
$$\equiv R_{1,d}(\theta) + R_{2,d}(\theta).$$

We first consider $R_{1,d}(\theta)$. On the event $A_d$, the signal is sparse and $\widehat{\theta}^*$ is the nonnegative garrote estimator $\widehat{\theta}_i^* = (1 - 2\log d/x_i^2)_+ x_i$ by the definition of $\widehat{\theta}^*$ in (8). Decomposing $R_{1,d}(\theta)$ further into two parts with either $\mu_d = d^{-1}\|\theta\|_2^2 \leq 3\gamma_d$ or $\mu_d > 3\gamma_d$ yields that $R_{1,d}(\theta) \leq R_F(\theta) I(\mu_d \leq 3\gamma_d) + r_{1,d}(\theta)$ where $R_F(\theta)$ is the risk of the nonnegative garrote estimator and $r_{1,d}(\theta) = d^{-1} E_\theta \{\|\widehat{\theta}^* - \theta\|_2^2 I(A_d)\} I(\mu_d > 3\gamma_d)$. The oracle inequality (3.10) in Cai (1991) with $L = 1$ and $\lambda = \lambda^F = 2\log d$ yields that

$$R_F(\theta) \leq d^{-1} \sum_i [(\theta_i^2 \wedge \lambda^F) + 4(\pi \log d)^{-1/2} d^{-1}]$$

(32)

$$\leq d^{-1} [\|\theta\|_2^2 + 4(\pi \log d)^{-1/2}].$$

Recall that $\mu_d = \|\theta\|_2^2/d$ and $\gamma_d = d^{-1/2} \log_2^{3/2} d$. It then follows from (32) that

(33) $$R_F(\theta) I(\mu_d \leq 3\gamma_d) \leq d^{-1}[3d^{1/2} \log_2^{3/2} d + 4(\pi \log d)^{-1/2}]$$
$$\leq cd^{-1/4}(\log d)^{5/2}.$$

Note that on the event $A_d$, $\|\widehat{\theta}^*\|_2^2 \leq \|x\|_2^2 \leq d + d\gamma_d$ and so

$$r_{1,d}(\theta) \leq 2d^{-1}(E\|\widehat{\theta}^*\|_2^2 + \|\theta\|_2^2) P(A_d) I(\mu_d > 3\gamma_d)$$
$$\leq 2(1 + 2\mu_d) P(A_d) I(\mu_d > 3\gamma_d) = o(d^{-1/2}),$$

where the last step follows from Lemma 2. Note that for any $\eta > 0$, Lemma 1 yields that

$$\widetilde{R}(\theta) - R_{block.oracle}(\theta) \leq C_\eta d^{\eta-1/2} \tag{34}$$

for some constant $C_\eta > 0$ and for all $\theta \in R^d$. Equations (27)–(29) yield

$$R_{2,d}(\theta) - \widetilde{R}(\theta) \leq \Delta_d \leq cd^{-1/4}(\log d)^{5/2}. \tag{35}$$

The proof of part (a) of the theorem is completed by putting together (33)–(35).

We now turn to part (b). It follows from part (a) that $R(\hat{\theta}^*, \theta) \leq R_{block.oracle}(\theta) + cd^{-1/4}(\log d)^{5/2}$. Note that $R_{block.oracle}(\theta) \leq r(d^{1/2} - 2, d^{1/2})$. Stein's unbiased risk estimate and Jensen's inequality yield that

$$r(d^{1/2} - 2, d^{1/2}) = d^{-1} \sum_b \left( d^{1/2} - (d^{1/2} - 2)^2 E \frac{1}{\|\underline{X}_b\|^2} \right)$$

$$\leq d^{-1} \sum_b \left( d^{1/2} - (d^{1/2} - 2)^2 \frac{1}{E\|\underline{X}_b\|^2} \right).$$

Note that $E\|\underline{X}_b\|^2 = \|\underline{\theta}_b\|^2 + d^{1/2}$. Hence $r(d^{1/2} - 2, d^{1/2}) \leq d^{-1} \sum_b (d^{1/2} + 1 - \frac{d}{\|\underline{\theta}_b\|^2 + d^{1/2}})$. The elementary inequality $(\sum_{i=1}^m a_i)(\sum_{i=1}^m a_i^{-1}) \geq m^2$, for $a_i > 0, 1 \leq i \leq m$ yields that

$$R_{block.oracle}(\theta) \leq r(d^{1/2} - 2, d^{1/2}) \leq d^{-1}\left( d + d^{1/2} - \frac{d^2}{\|\theta\|_2^2 + d} \right)$$

$$= \frac{\|\theta\|_2^2}{\|\theta\|_2^2 + d} + d^{-1/2}$$

and part (b) then follows. We now consider part (c). Note first that $R_{1,d}(\theta) \leq R_F(\theta)$. On the other hand, $R_{2,d}(\theta) = d^{-1} E\{\|\hat{\theta}^* - \theta\|_2^2 I(A_d^c)\} \leq Cd^{-1}(E\|\hat{\theta}^* - \theta\|_2^4)^{1/2} P^{1/2}(A_d^c)$. To complete the proof of (12), it then suffices to show that under the assumption $\mu_d \leq \frac{1}{3}\gamma_d$, $E\|\hat{\theta}^* - \theta\|_2^4$ is bounded by a polynomial of $d$ and $P(A_d^c)$ decays faster than any polynomial of $d^{-1}$. Note that in this case $\|\theta\|_2^2 = d\mu_d \leq d^{1/2} \log_2^{3/2} d$. Since $\|\hat{\theta}^*\|_2^2 \leq \|x\|_2^2$ and $x_i = \theta_i + z_i$,

$$E\|\hat{\theta}^* - \theta\|_2^4 \leq E(2\|\hat{\theta}^*\|_2^2 + 2\|\theta\|_2^2)^2$$

$$\leq E(2\|x\|_2^2 + 2\|\theta\|_2^2)^2 \leq E(4\|\theta\|_2^2 + 2\|z\|_2^2)^2$$

$$\leq 32\|\theta\|_2^4 + 8E\|z\|_2^4 \leq 32d\log_2^2 d + 16d + 8d^2.$$

On the other hand, it follows from Hoeffding's inequality and Mill's inequality that

$$P(A_d^c) = P\left( d^{-1} \sum (z_i^2 + 2z_i\theta_i + \theta_i^2 - 1) > d^{-1/2} \log_2^{3/2} d \right)$$

$$\leq P\Big(d^{-1}\sum(z_i^2-1)>\tfrac{1}{3}d^{-1/2}\log_2^{3/2}d\Big)$$

$$+ P\Big(d^{-1}\sum 2z_i\theta_i>\tfrac{1}{3}d^{-1/2}\log_2^{3/2}d\Big)$$

$$\leq 2\exp(-C\log_2^3 d)+\tfrac{1}{2}\exp(-C\log_2^3 d/\mu_d),$$

which decays faster than any polynomial of $d^{-1}$.

6.3. *Proof of Theorem 2.* Note that $x/(x+d)$ is increasing in $x$ and for $p\geq 2$, $\sup_{\theta\in\Theta_p(\tau)}\|\theta\|_2=d^{1/2-1/p}\tau$. Hence

$$\sup_{\theta\in\Theta_p(\tau)} R(\hat{\theta}^*,\theta)\leq \sup_{\theta\in\Theta_p(\tau)}\frac{\|\theta\|_2^2}{\|\theta\|_2^2+d}+cd^{-1/4}(\log d)^{5/2}$$

$$\leq \frac{\sup_{\theta\in\Theta_p(\tau)}\|\theta\|_2^2}{\sup_{\theta\in\Theta_p(\tau)}\|\theta\|_2^2+d}+cd^{-1/4}(\log d)^{5/2}$$

$$=\frac{d^{1-2/p}\tau^2}{d^{1-2/p}\tau^2+d}+cd^{-1/4}(\log d)^{5/2}$$

$$=\frac{\tau^2}{\tau^2+d^{2/p}}+cd^{-1/4}(\log d)^{5/2}.$$

Now consider part (b). It follows from Proposition 3 that there is a constant $c_p$ depending on $p$ such that $R(\theta)\leq \inf_{\lambda\geq 0}r(\lambda,1)\leq cd^{-1}\tau^p(\log(d\tau^{-p}))^{(2-p)/2}$, since $\Theta_p(\tau)=\{\theta\in\mathbb{R}^d:\|\theta\|_p^p/d\leq\tau^p/d\}$. Part (b) now follows directly from (11) in Theorem 1.

For part (c) it is easy to check that for $0<p\leq 2$, $\|\theta\|_2^2\leq\|\theta\|_p^2\leq\tau^2<\tfrac{1}{3}d^{1/2}\log_2^{3/2}d$ and so $\mu_d\leq\tfrac{1}{3}\gamma_d$. It then follows from Theorem 1 and (39) that

$$R(\hat{\theta}^*,\theta)\leq R_F(\theta)+cd^{-1}(\log d)^{-1/2}$$

$$\leq \frac{1}{d}(\|\theta\|_2^2+8(2\log d)^{-1/2})+cd^{-1}(\log d)^{-1/2}$$

$$\leq \frac{1}{d}\tau^2+cd^{-1}(\log d)^{-1/2}.$$

6.4. *Proof of Theorem 3.* Again we set $\sigma=1$. Note that the empirical wavelet coefficients $\tilde{y}_{j,k}$ can be written as $\tilde{y}_{j,k}=\tilde{\theta}_{j,k}+n^{-1/2}z_{j,k}$, where $\tilde{\theta}_{j,k}=E\tilde{y}_{j,k}$ are the DWT of the sampled function $\{n^{-1/2}f(\tfrac{i}{n})\}$ and $z_{j,k}\overset{\text{i.i.d.}}{\sim}N(0,1)$. To more conveniently use Theorem 1, we multiply both sides by $n^{1/2}$ and get

$$(36)\qquad y'_{j,k}=\theta'_{j,k}+z_{j,k},\qquad j\geq j_0,\ k=1,2,\ldots,2^j,$$

where $y'_{j,k} = n^{1/2}\tilde{y}_{j,k}$ and $\theta'_{j,k} = n^{1/2}\tilde{\theta}_{j,k}$. Let $\tilde{f} = \sum_{k=1}^{2^J} n^{-1/2} f(\frac{i}{n})\phi_{J,k}(t)$, where $n = 2^J$. Note that $\sup_{f \in B^{\alpha}_{p,q}(M)} \|f - \tilde{f}\|_2^2 = o(n^{-2\alpha/(1+2\alpha)})$ from Lemma 4. To establish Theorem 3, by the Cauchy–Schwarz inequality as in Remark 1, it suffices to show that

$$\sup_{f \in B^{\alpha}_{p,q}(M)} E_f \|\hat{f}^* - \tilde{f}\|_2^2 \leq R_T^*(B^{\alpha}_{p,q}(M))(1 + o(1)).$$

Fix $0 < \varepsilon_0 < 1/(1 + 2\alpha)$ and let $J_0$ be the largest integer satisfying $2^{J_0} \leq n^{\varepsilon_0}$. Write

$$n^{-1} E_{\theta'} \|\hat{\theta}'^* - \theta'\|_2^2 = \left( \sum_{j \leq J_0} \sum_k + \sum_{J_0 \leq j < J_1} \sum_k + \sum_{j \geq J_1} \sum_k \right) n^{-1} E_{\theta'} \|\hat{\theta}'^* - \theta'\|_2^2$$

$$= S_1 + S_2 + S_3,$$

where $J_1 > J_0$ is to be chosen later. The terms $S_1$ and $S_3$ are identical in all block thresholding procedures. We thus only focus on the term $S_2$. Since $0 < \varepsilon_0 < 1/(1 + 2\alpha)$, $n^{\varepsilon_0} n^{-1} \log n = o(n^{-2\alpha/(1+2\alpha)})$. It then follows from (26) in Lemma 3 that $S_1 \leq C n^{\varepsilon_0} n^{-1} \log n = o(n^{-(2\alpha)/(1+2\alpha)})$ which is negligible relative to the minimax risk. On the other hand, it follows from (11) in Theorem 1 that

$$S_2 \leq \sum_{J_0 \leq j < J_1} n^{-1} 2^j R(\underline{\theta}'_j) + \sum_{J_0 \leq j < J_1} n^{-1} 2^j R_F(\underline{\theta}'_j) I(\mu'_{2^j} \leq 3\gamma_{2^j})$$

$$+ \sum_{J_0 \leq j < J_1} n^{-1} c 2^{3j/4} j^{5/2}$$

$$= S_{21} + S_{22} + S_{23},$$

where $\mu'_{2^j} = 2^{-j} \|\underline{\theta}'_j\|_2^2 = 2^{-j} n \|\underline{\tilde{\theta}}_j\|_2^2$ and $\gamma_{2^j} = 2^{-j/2} j^{3/2}$. It follows from Remark 1 that

$$(37) \qquad S_{21} = \sum_{J_0 \leq j < J_1} n^{-1} 2^j R(\underline{\theta}'_j) \leq (1 + o(1)) R_T^*(B^{\alpha}_{p,q}(M)).$$

We shall see that both $S_{22}$ and $S_{23}$ are negligible relative to the minimax risk. Note that

$$(38) \qquad S_{23} = \sum_{J_0 \leq j < J_1} n^{-1} c 2^{3j/4} j^{5/2} \leq C n^{-1} 2^{3J_1/4} J_1^{5/2}.$$

The oracle inequality (3.10) in Cai (1991) with $L = 1$ and $\lambda = \lambda^F = 2\log(2^j)$ yields that

$$2^j R_F(\underline{\theta}'_j) \leq \sum_{k=1}^{2^j} (n \tilde{\theta}_{j,k} \wedge \lambda^F) + 8(2\log 2^j)^{-1/2} 2^{-j}$$

$$(39)$$

$$\leq n \|\underline{\tilde{\theta}}_j\|_2^2 + 8(2\log 2^j)^{-1/2} 2^{-j}.$$

Recall that $\mu'_{2^j} = 2^{-j} n \|\tilde{\underline{\theta}}_j\|_2^2$ and $\gamma_{2^j} = 2^{-j/2} j^{3/2}$. It then follows from (39) that

$$
\begin{aligned}
(40) \quad S_{22} &= \sum_{J_0 \le j < J_1} n^{-1} 2^j R_F(\underline{\theta}'_j) I(\mu'_{2^j} \le 3\gamma_{2^j}) \\
&\le \sum_{J_0 \le j < J_1} (3n^{-1} 2^{j/2} j^{3/2} + 8n^{-1}(2\log 2^j)^{-1/2} 2^{-j}) \le Cn^{-1} 2^{J_1/2} J_1^{3/2}.
\end{aligned}
$$

Hence if $J_1$ satisfies $2^{J_1} = n^\gamma$ with some $\gamma < \frac{4}{3(1+2\alpha)}$, then (38) and (40) yield

$$
(41) \qquad S_{22} + S_{23} = o(n^{-(2\alpha)/(1+2\alpha)}).
$$

We now turn to the term $S_3$. It is easy to check that for $\theta \in B^\alpha_{p,q}(M)$, $\|\underline{\theta}_j\|_2^2 \le M^2 2^{-2\alpha' j}$ where $\alpha' = \alpha - (\frac{1}{p} - \frac{1}{2})_+ > 0$. Note that if $J_1$ satisfies $2^{J_1} = n^\gamma$ for some $\gamma > \frac{1}{1/2+2\alpha'}$, then for all sufficiently large $n$ and all $j \ge J_1$, $2^{-j} n \|\underline{\theta}_j\|_2^2 \le \frac{1}{4}\gamma_{2^j}$ where $\gamma_{2^j} = 2^{-j/2} j^{3/2}$ which implies $\mu'_{2^j} \le \frac{1}{3}\gamma_{2^j}$ for $j \ge J_1$ and $\gamma > 2(1-2\beta)$, since

$$
(42) \quad \mu'_{2^j} - 2^{-j} n \|\underline{\theta}_j\|_2^2 \le 2^{-j} n \|\tilde{\underline{\theta}}_j - \underline{\theta}_j\|_2^2 \le C2^{-j} n^{1-2\beta} = o(2^{-j/2} j^{3/2}).
$$

It thus follows from (12) and (39) that

$$
\begin{aligned}
(43) \quad S_3 &= \sum_{j \ge J_1} \sum_k n^{-1} E_{\theta'}(\hat{\theta}'^*_{j,k} - \theta'_{j,k})^2 \\
&\le \sum_{j \ge J_1} (n^{-1} 2^j R_F(\underline{\theta}'_j) + cn^{-1} j^{-j/2}) \\
&\le \sum_{j \ge J_1} \|\tilde{\underline{\theta}}_j\|_2^2 + Cn^{-1} \\
&\le \sum_{j \ge J_1} \|\underline{\theta}_j\|_2^2 + Cn^{-1} + \|\tilde\theta - \theta\|_2^2 \le C2^{-2\alpha' J_1} + Cn^{-1} \\
&= o(n^{-(2\alpha)/(1+2\alpha)}),
\end{aligned}
$$

when $\gamma > \frac{\alpha}{(1+2\alpha)\alpha'}$ and $\beta > \frac{\alpha}{1+2\alpha}$. Equations (41)–(43) hold by choosing $\gamma$ satisfying

$$
\max\left\{ \frac{\alpha}{(1+2\alpha)\alpha'}, 2(1-2(\alpha-1/p)), \frac{1}{1/2+2\alpha'} \right\} < \gamma < \frac{4}{(1+2\alpha)3},
$$

which is possible for $\alpha > 4(\frac{1}{p} - \frac{1}{2})_+ + \frac{1}{2}$ with $\frac{2\alpha^2 - 1/6}{1+2\alpha} > 1/p$. This completes the proof.

6.5. *Proof of Theorem 4.* Define the minimax linear risk by

$$R_L^*(B_{p,q}^\alpha(M)) = \inf_{\widehat{\theta} \text{ linear}} \sup_{\theta \in B_{p,q}^\alpha(M)} E\|\widehat{\theta} - \theta\|_2^2.$$

It follows from Donoho, Liu and MacGibbon ([1990](#)) and Remark [1](#) that

$$R_L^*(B_{p,q}^\alpha(M)) = (1 + o(1)) \sup_{\theta \in B_{p,q}^\alpha(M)} \sum_{j=j_0}^\infty \sum_k \left( \frac{\theta_{j,k}^2/n}{\theta_{j,k}^2 + 1/n} \right),$$

$R_L^*(B_{p,q}^\alpha(M)) = R^*(B_{p,q}^\alpha(M))(1 + o(1))$ for $p = q = 2$, and $R_L^*(B_{p,q}^\alpha(M)) \leq 1.25 R^*(B_{p,q}^\alpha(M))(1 + o(1))$ for all $p \geq 2$ and $q \geq 2$. It thus suffices to show that SureBlock asymptotically attains the minimax linear risk for $\alpha > 0$, $p \geq 2$ and $q \geq 2$. Since $\|\theta - \tilde{\theta}\|_2^2 = o(n^{-2\beta})$ with $\beta > \alpha/(1 + 2\alpha)$, we need only to show $\sup_{\theta \in B_{p,q}^\alpha(M)} E_\theta \|\hat{\theta}^* - \tilde{\theta}\|_2^2 \leq R_L^*(B_{p,q}^\alpha(M))(1 + o(1))$ similar to the arguments in Remark [1](#).

Recall in the proof of Theorem [3](#) it is shown that $E_\theta \|\hat{\theta}^* - \tilde{\theta}\|_2^2 \leq S_1 + S_{21} + S_{22} + S_{23} + S_3$, where $S_1 + S_{22} + S_{23} + S_3 = o(n^{-2\alpha/(2\alpha+1)})$ and $S_{21} = \sum_{J_0 \leq j < J_1} n^{-1} 2^j R(\underline{\theta}'_j)$ with $J_0$ and $J_1$ chosen as in the proof of Theorem [3](#). Since the minimax risk $R^*(B_{p,q}^\alpha(M)) \asymp n^{-2\alpha/(2\alpha+1)}$, this implies that $S_{21}$ is the dominating term in the maximum risk of SureBlock. It follows from the definition of $R(\underline{\theta}'_j)$ given in ([10](#)) that $n^{-1} 2^j R(\underline{\theta}'_j) \leq n^{-1} \sum_b E_{\underline{\theta}'_b} \|\hat{\underline{\theta}}'_b(L_j - 2, L_j) - \underline{\theta}'_b\|_2^2$, where the RHS is the risk of the blockwise James–Stein estimator with any fixed block size $1 \leq L_j \leq 2^{j/2}$ and a fixed threshold level $L_j - 2$. Stein's unbiased risk estimate [see, e.g., Johnstone ([2002](#)), Chapter 9.2] yields that $n^{-1} \sum_b E_{\underline{\theta}'_b} \|\hat{\underline{\theta}}'_b(L_j - 2, L_j) - \underline{\theta}'_b\|_2^2 \leq \sum_b \left( \frac{\|\hat{\underline{\theta}}_b\|_2^2 L_j/n}{\|\hat{\underline{\theta}}_b\|_2^2 + L_j/n} + \frac{2}{n} \right)$. Hence the maximum risk of SureBlock satisfies

$$\sup_{\theta \in B_{p,q}^\alpha(M)} E_\theta \|\hat{\theta}^* - \tilde{\theta}\|_2^2$$

$$\leq \sup_{\theta \in B_{p,q}^\alpha(M)} \sum_{J_0 \leq j < J_1} \sum_b \left( \frac{\|\tilde{\underline{\theta}}_b\|_2^2 L_j/n}{\|\hat{\underline{\theta}}_b\|_2^2 + L_j/n} + \frac{2}{n} \right) \cdot (1 + o(1))$$

$$\leq \sup_{\theta \in B_{p,q}^\alpha(M)} \sum_{J_0 \leq j < J_1} \sum_b \left( \frac{\|\underline{\theta}_b\|^2 L_j/n}{\|\underline{\theta}_b\|^2 + L_j/n} + \frac{2}{n} \right) \cdot (1 + o(1))$$

$$= \sup_{\theta \in B_{p,q}^\alpha(M)} \left( \sum_{J_0 \leq j < J_1} \sum_b \frac{\|\underline{\theta}_b\|_2^2 L_j/n}{\|\underline{\theta}_b\|_2^2 + L_j/n} + 2 \sum_{J_0 \leq j < J_1} \frac{2^j}{n L_j} \right) \cdot (1 + o(1)),$$

where the second inequality follows from a similar argument as in Remark [1](#). Note that in the proof of Theorem [1](#), $J_1$ satisfies $2^{J_1} = n^\gamma$ with $\gamma < \frac{4}{(1+2\alpha)3}$.

Hence if $L_j$ satisfies $2^{j\rho} \le L_j \le 2^{j/2}$ for some $\rho > \frac{1}{4}$, then $\sum_{j \le J_1} \frac{2^j}{nL_j} \le \frac{1}{n} 2 \cdot 2^{J_1 3/4} = o(n^{-(2\alpha)/(1+2\alpha)})$ and hence

(44)
$$\sup_{\theta \in B_{p,q}^\alpha(M)} E_\theta \|\hat{\theta}^* - \theta\|_2^2$$

$$\le \sup_{\theta \in B_{p,q}^\alpha(M)} \sum_{J_0 \le j < J_1} \sum_b \left( \frac{\|\underline{\theta}_b\|^2 L_j/n}{\|\underline{\theta}_b\|^2 + L_j/n} \right) \cdot (1 + o(1)).$$

Note that $\sum_{b=1}^{2^j/L_j} \frac{\|\underline{\theta}_b\|^2 L_j/n}{\|\underline{\theta}_b\|^2 + L_j/n} = \frac{2^j}{n} - (\frac{L_j}{n})^2 \sum_{b=1}^{2^j/L_j} \frac{1}{\|\underline{\theta}_b\|^2 + L_j/n}$. Then the simple inequality $(\sum_{i=1}^m a_i)(\sum_{i=1}^m a_i^{-1}) \ge m^2$, for $a_i > 0, 1 \le i \le m$ yields that

(45)
$$\sum_{J_0 \le j < J_1} \sum_b \left( \frac{\|\underline{\theta}_b\|^2 L_j/n}{\|\underline{\theta}_b\|^2 + L_j/n} \right) \le \frac{2^j}{n} - \left( \frac{L_j}{n} \right)^2 \left( \frac{2^j}{L_j} \right)^2 \frac{1}{\sum_{b=1}^{2^j/L_j} \left( \|\underline{\theta}_b\|^2 + \frac{L_j}{n} \right)}$$

$$= \frac{2^j/n \sum_{b=1}^{2^j/L_j} \|\underline{\theta}_b\|^2}{\sum_{b=1}^{2^j/L_j} \|\underline{\theta}_b\|^2 + 2^j/n}$$

$$= \frac{2^j/n \sum_k |\theta_{j,k}|^2}{\sum_k |\theta_{j,k}|^2 + 2^j/n}.$$

Theorem 2 in Cai, Low and Zhao (2000) shows that

(46)
$$\sup_{\theta \in B_{p,q}^\alpha(M)} \sum_{J_0 \le j < J_1} \frac{2^j/n \sum_k |\theta_{j,k}|^2}{\sum_k |\theta_{j,k}|^2 + 2^j/n} = R_L^*(B_{p,q}^\alpha(M))(1 + o(1)).$$

The proof is complete by combining (44)–(46).

6.6. *Proof of Theorem 5.* Set $\rho_g(\eta) \triangleq \inf_\lambda \sup_{\mathcal{F}_p(\eta)} E_F r_g(\mu)$ where $\mathcal{F}_p(\eta)$ and $r_g(\mu)$ are given as in Proposition 3. Proposition 3 implies that $\rho_g(\eta) \le \overline{r}(\delta_\lambda^g, \eta) \le 2\eta^p (2 \log \eta^{-p})^{(2-p)/2} (1 + o(1))$ as $\eta \to 0$. For $p \in (0, 2)$, Theorem 15 of Donoho and Johnstone (1994b) shows the univariate Bayes minimax risk satisfies $\rho(\eta) \triangleq \inf_\delta \sup_{\mathcal{F}_p(\eta)} E_F E_\mu (\delta(x) - \mu)^2 = \eta^p (2 \log \eta^{-p})^{(2-p)/2} (1 + o(1))$ as $\eta \to 0$. Note that $\rho_g(\eta)/\rho(\eta)$ is bounded as $\eta \to 0$ and $\rho_g(\eta)/\rho(\eta) \to 1$ as $\eta \to \infty$. Both $\rho_g(\eta)$ and $\rho(\eta)$ are continuous on $(0, \infty)$, so $G(p) = \sup_\eta \frac{\rho_g(\eta)}{\rho(\eta)} < \infty$, for $p \in (0, 2)$. Theorems 4 and 5 in Section 4 of Donoho and Johnstone (1998) derived the asymptotic minimaxity over Besov bodies from the univariate Bayes minimax estimators. It then follows from an analogous argument of Section 5.3 in Donoho and Johnstone (1998) that

$$R_T^*(B_{p,q}^\alpha(M)) \le \inf_{\lambda_j} \sup_{\theta \in B_{p,q}^\alpha(M)} E \sum_{j=j_0}^\infty \|\hat{\underline{\theta}}_j(\lambda_j, 1) - \underline{\theta}_j\|^2$$

$$\le G(p \wedge q) \cdot R^*(B_{p,q}^\alpha(M))(1 + o(1)).$$

## APPENDIX: TEST FUNCTIONS

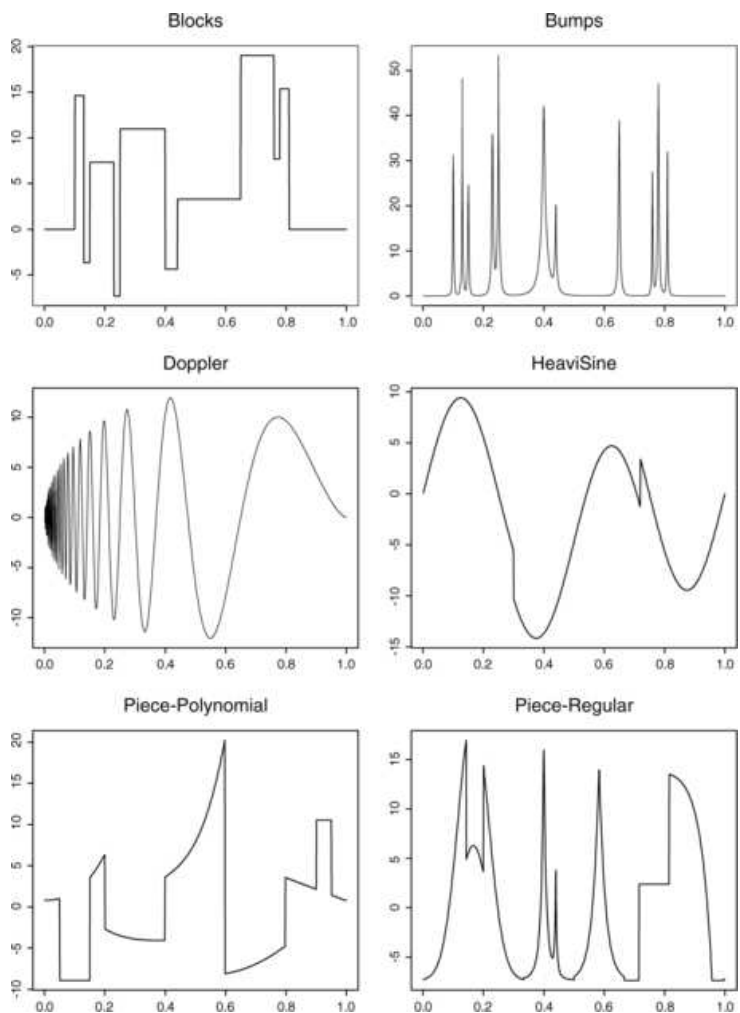

Fig. A.1.  *Test functions.*

DEPARTMENT OF STATISTICS
THE WHARTON SCHOOL
UNIVERSITY OF PENNSYLVANIA
PHILADELPHIA, PENNSYLVANIA 19104
USA
E-MAIL: tcai@wharton.upenn.edu

DEPARTMENT OF STATISTICS
YALE UNIVERSITY
NEW HAVEN, CONNECTICUT 06511
USA
E-MAIL: huibin.zhou@yale.edu